\documentclass[11pt]{article}
 \usepackage{indentfirst}  
\usepackage[a4paper,margin=1in]{geometry}
\usepackage[T1]{fontenc}
\usepackage[utf8]{inputenc}
\usepackage{lmodern}
\usepackage{amsmath,amssymb,amsfonts,amsthm,mathrsfs}
\usepackage{mathtools}
\usepackage{enumitem}
\usepackage{microtype}
\usepackage{hyperref}
\usepackage[nameinlink,capitalize]{cleveref}

\hypersetup{
  colorlinks=true,
  linkcolor=blue,
  citecolor=blue,
  urlcolor=blue
}

\newtheorem{theorem}{Theorem}[section]
\newtheorem{lemma}[theorem]{Lemma}

\newtheorem{corollary}[theorem]{Corollary}
\theoremstyle{definition}

\newtheorem{remark}[theorem]{Remark}

\newcommand{\OO}{\mathcal O}

\newcommand{\SSp}{\mathbb S}
\newcommand{\pt}{\operatorname{pt}}
\newcommand{\sat}{\operatorname{sat}}
\newcommand{\cl}{\operatorname{cl}}
\newcommand{\X}{\operatorname{X}}
\newcommand{\Q}{\operatorname{Q}}
\newcommand{\Sig}{\Sigma}
\newcommand{\BoxU}{\mathord\Box}

\begin{document}
\renewcommand{\thefootnote}{\alph{footnote}}

\hspace*{7mm} \setcounter{footnote}{-1} \footnote{\textit{$*$ Corresponding author.\\
E-mail
address}: liubei@hebtu.edu.cn(B.Liu),  liqingguoli@aliyun.com(Q.Li), lcx19882012@163.com(C.Lu)}
\begin{center}

{\bf \Large Function spaces with compact-open topology  need not be sober}

\vskip 0.3cm \centerline{ Bei Liu$^1$, Qingguo Li$^{2*}$, Chongxia Lu$^3$} {\small \textit{
$^1$ School of Mathematical Sciences, Hebei Normal University, Shijiazhuang  050024, China\\
$^2$  School of Mathematics, Hunan University, Changsha 410082, China\\
$^3$ School  of Mathematics and Statistics, Ningbo University, Ningbo 315211, China}\\}

 \vskip0.3cm
\date{}
 {\begin{minipage}{135mm}\small {\bf Abstract:}
We exhibit sober spaces $X$ and $Y$ such that the function space $[X \to Y]_K$ of continuous maps, equipped with the compact-open topology, is not sober; the codomain can be chosen to be the Sierpiński space $S$. Our construction builds on Isbell's complete lattice $M$, whose Scott space $\Sigma M$ is non-sober. Let $H = Q(M)$ be the set of nonempty compact saturated subsets of $\Sigma M$, ordered by reverse inclusion. Since $H$ is well known to be a spatial frame, the space $X = \mathrm{pt}(H)$ is sober with $\mathcal{O}(X) \cong H$. The key technical result is that sobriety of the space $H$ would entail sobriety of the original Scott space $\Sigma M$. This contradicts Isbell's example, and the desired counterexample follows. 

\vspace*{12pt}{\bf Keywords:}
Compact-open topology; sober; Isbell's complete lattice; Scott open filter

\end{minipage}}
\end{center}

\vskip0.2in\section{Introduction}

Sobriety is a central separation property in non-Hausdorff topology and domain theory.    In locale theory, sober spaces are the point-set representatives of spatial frames, and in domain theory,  sobriety is one of the fundamental links between order-theoretic and topological semantics, for example, \cite{GierzEtAl2003,GoubaultLarrecq2013}.

Function spaces provide a natural test of the stability of topological properties.  If $X$ and $Y$ are topological spaces, the set $\mathbb{C}(X,Y)$ of continuous maps can be equipped with several familiar topologies, including the pointwise convergence topology, Isbell topology, and compact-open topology.  For weak variants of sobriety, a number of preservation and reflection questions for function spaces have been studied (see \cite{Ershov,Jia,LiuLiHo2021,xubao,ZhangBaoXu2022}).  The compact-open topology is particularly natural because it is generated by conditions on compact subsets of the domain.

The purpose of this paper is to show that sobriety is not preserved by the compact-open function space, even when both the domain and codomain are sober.  More precisely, we prove the following statement:

There exists a sober space $X$ such that, for the Sierpi\'nski space $\SSp$, the compact-open function space \([X\rightarrow\SSp]_K\)
is not sober.  In particular, there exist sober spaces $X$ and $Y$ for which $[X\rightarrow Y]_K$ is non-sober. Thus, giving an affirmative answer to the problem posed in \cite{xubao}.

The construction uses three ingredients:
\begin{enumerate}[label=\textup{(\roman*)}]
\item Isbell constructed a complete lattice $M$ for which the Scott space $\Sig M$ is not sober \cite{Isbell1982}; see also the earlier example of Johnstone \cite{Johnstone1981}.
\item If $M$ is a complete lattice, the family $\Q(M)$ of nonempty compact saturated subsets of $\Sig M$, ordered by the reverse inclusion, is a complete Heyting algebra; in fact,  it is spatial.  This construction is treated explicitly in \cite{XuXiZhao2021}.  The proof relies on well-filteredness of Scott spaces of complete lattices \cite{XiLawson2017} and coherence results for such dcpos \cite{JiaJungLi2016}.
\item Every spatial frame $H$ is isomorphic to the open set frame of its point space $\pt(H)$, and that point space is sober.
\end{enumerate}

The main additional point is   to show that sobriety of the resulting topology on $\Q(M)$ would imply sobriety of $\Sig M$.  We prove this transfer statement directly.  The proof recovers a least element from a hypothetical generic compact saturated set and then shows that any nonempty irreducible Scott closed set in $M$ must be principal.  Applied to Isbell's lattice, this gives the required contradiction.

The paper is organized as follows.  In \cref{sec:preliminaries},  we collect the order-topological notions and useful results needed in the proof. In \cref{sec:frame},   we discuss the spatial   frame $\Q(M)$, and pass from spatial frames to sober spaces.   The main theorem and an explicit witness to non-sobriety are given in \cref{sec:main-result}.

\section{Preliminaries}\label{sec:preliminaries}

In this section, we briefly recall some basic concepts and notations to be used in this paper; more details can be found in \cite{gelun,GierzEtAl2003,GoubaultLarrecq2013,Johnstone1982}.

For a poset $P$ and $A \subseteq P$, set
$\uparrow A = \{x \in P : a \leq x \text{ for some } a \in A\}$ (dually $\downarrow A = \{x \in P : a \geq x \text{ for some } a \in A\}$). For
$p\in P$, we abbreviate $\uparrow\{p\}$ to be $\uparrow p$. 
A subset $A$ is called an upper set (resp., a lower set) if $A=\uparrow A$ (resp., $A=\downarrow A$).
If the set of
upper bounds of $A$ has a least element, that element is the supremum of
$A$, denoted by $\bigvee A$ or $\sup A$.  The greatest lower bound of $A$ is also called the infimum of $A$ and is denoted by $\bigwedge A$. If $\bigvee A$ and $\bigwedge A$ exist for all $A \subseteq P$, then $P$ is called a complete lattice. 
Let $\bigvee P=1_P$ and $\bigwedge P=0_P$.
For a poset $P$, a subset $D$ of $P$ is \emph{directed} provided that it is non-empty, and every finite subset of $D$ has an upper bound in $D$.
A subset $U\subseteq P$ is  \emph{Scott open}  if (1) $U$ is an upper set,  and
(2)  whenever $D\subseteq P$ is directed, $\bigvee D$ exists, and $\bigvee D\in U$, one has $D\cap U\neq\varnothing$.
The Scott topology is denoted by $\sigma(P)$ and the Scott space by $\Sig P=(P,\sigma(P))$.

Given a $T_0$ space $X$, we can define a partial order $x\leq_{\X} y$, called  the specialization order, which is  defined by
\(x\leq_{\X}  y\) if and only if is \(x\in\cl_X\{y\}\). Let $\OO(X)$ be the set of all open subsets of $X$. 
A subset $A\subseteq X$ is \emph{saturated} if it is an intersection of open subsets of $X$.  
A nonempty subset $A\subseteq X$   is  called \emph{irreducible} if, whenever $U$ and $V$ are open and both meet $A$, then $U\cap V$ also meet $A$.   
A $T_0$ space $X$ is \emph{sober} if every nonempty irreducible closed subset $A$ has a unique generic point, i.e., a unique $x\in X$ such that
\(A=\cl_X\{x\}\). 
We denote by $\Q(X)$ the poset of non-empty compact saturated subsets of $X$ with the reverse inclusion order, i.e., for $K_1,K_2\in\Q(X)$, $K_1\sqsubseteq K_2$ iff $K_2\subseteq K_1$. 
The space $X$ is called \emph{well-filtered} if for any filtered family $\mathcal{K}\subseteq\Q(X)$ and $U\in\OO(X)$, $\bigcap\mathcal{K}\subseteq U$ implies $K\subseteq U$ for some $K\in\mathcal{K}$. For $U\in\OO(X)$, let $\BoxU U=\{K\in\Q(X):K\subseteq U\}$. The upper Vietoris topology on $\Q(X)$ is the topology that has $\{\BoxU U:U\in\OO(X)\}$ as a base and the space is denoted by $P_S(X)$.
For $G\subseteq X$, let $\Phi(G)=\{U\in\OO(X):G\subseteq U\}$. Let \(\phi(\OO(X))\) be the topology on $\OO(X)$ generated by $\{\Phi(C):  C \ \text{is a compact subset of}\ X\}$. It is easy to verify that the family is a base of  $\{\Phi(K):K\in\Q(X)\cup\emptyset\}$.
For a complete lattice $L$, let $\Q(L)$ be the collection of all non-empty saturated compact saturated subsets of   $\Sigma L$.

A  \emph{frame} is a complete lattice $H$ satisfying the infinite distributive law: 
\[h\wedge\bigvee_{i\in I}h_i=\bigvee_{i\in I}(h\wedge h_i),
\] 
 for every $h\in H$ and every family $\{h_i\}_{i\in I}\subseteq H$.  Such a   frame is also called a complete Heyting algebra. 
  It is well-known that a complete lattice is a \emph{spatial} frame if and only if it is isomorphic to the lattice of all open subsets of some topological space (cf.\cite{Johnstone1982}).
 Let $\mathbf 2=\{0,1\}$ be the two-element frame.  A  \emph{point}  of a frame $H$ is a frame homomorphism \(p:H\rightarrow\mathbf 2\), 
i.e., a map preserving arbitrary joins and finite meets.  Write
\[
   \pt(H)=\{p:H\to\mathbf 2:p\text{ is a frame homomorphism}\}.
\]
For $h\in H$, define
\[\widehat h=\{p\in\pt(H):p(h)=1\}.\]

\begin{lemma}
\label{lem:frame-spectrum}
(\cite{GierzEtAl2003})
Let $H$ be a spatial frame and $X=\pt(H)$, then
\begin{enumerate}
[label=\textup{(\roman*)}]
  \item the family $\{\widehat h:h\in H\}$ is a topology on $X$.
  \item The mapping \[\eta:H\rightarrow\OO(X),
   \qquad
   \eta(h)=\widehat h,\]
is a frame isomorphism, where $\OO(X)=\{\widehat h:h\in H\}$.
  \item The space $(X,\OO(X))$ is sober.
\end{enumerate}
\end{lemma}

\begin{lemma}
(\cite{GierzEtAl2003})
 \label{scottof}
Let $X$ be a $T_0$ space. If $K\subseteq X$ is compact, then
\[
\Phi(K)
=
\{U\in\mathcal{O}(X): K\subseteq U\}
\]
is a Scott open filter in $\mathcal{O}(X)$.
\end{lemma}

 \begin{remark}
 \label{sof}
 The set $K$ need not be saturated, so define
\[
   \sat(K)
   =
   \bigcap\{U\in\OO(X):K\subseteq U\}.
\]
Then it is easy to prove that
\[
   \Phi(K)=\Phi(\sat(K)),
\]
and  $\sat(K)$ is compact and saturated.
\end{remark}

\begin{lemma}
\label{lem:smyth-cosar-serscott}
(\cite{GierzEtAl2003})
Let $X$ be a well-filtered  space. Then the upper Vietoris topology is coarser than the Scott topology on $\Q(X)$.
\end{lemma}

\begin{lemma}
\label{lem:xilawson}
(\cite{XiLawson2017})
Let $L$ be a complete lattice or a bounded complete dcpo. Then $L $ equipped with the Scott
topology is well-filtered.
\end{lemma}

\begin{lemma}
\label{Ishey}
(\cite{XuXiZhao2021})
For any complete lattice L, $\Q(L)$ is a spatial frame.
\end{lemma}

Let $H$ be a frame.  We denote by $\tau_{\mathrm f}$ the topology on the underlying set of $H$ generated by all  Scott open proper   filters of $H$.
According to \cref{sof} and  \cref{lem:frame-spectrum} (ii), we have the following  corollary. 

\begin{corollary}\label{lem:filter-co}
Let $H$ be a spatial frame, let
$X=\pt(H)$
be the sober space of \cref{lem:frame-spectrum}, and $\tau_{\mathrm f}$   the   topology on $H$.  Then
\[
   (H,\tau_{\mathrm f})\cong (\OO(X),\phi(\OO(X))).
\]
\end{corollary}

\section{Some results related to complete lattice $M$ and $\Q(M)$}
\label{sec:frame}

We now specialize  $H=\Q(M)$
for a complete lattice $M$, always ordered by the reverse inclusion and equipped with the  topology $\tau_{\mathrm f}$.

\begin{lemma}
\label{lem:smyth-filter-properties}
Let $M$ be a complete lattice,  $H=\Q(M)$, 
and   $H$ endowed with the   topology $\tau_{\mathrm f}$.  Then the following hold:
\begin{enumerate}[label=\textup{(\roman*)}]
\item For every Scott open $U\subseteq M$, the set $\BoxU U$ is $\tau_{\mathrm f}$-open.
\item  The specialization order of $(H,\tau_{\mathrm f})$ is the reverse inclusion:
\[
   K\leq_{\tau_f} L
   \quad\Longleftrightarrow\quad
   K\supseteq L.
\]
\item The canonical map
\[
   \xi:\Sig M\rightarrow(H,\tau_{\mathrm f}),
   \qquad
   \xi(x)=\mathord\uparrow x,
\]
is continuous.
\end{enumerate}
\end{lemma}

\begin{proof}
(i).
If $U=\varnothing$, then $\BoxU U=\varnothing$.  If $U=M$, we get  $\BoxU U=H$.  Both are open.  

Therefore, we assume that $U$ is nonempty and proper.  We first prove that $\BoxU U$ is a proper filter of $H$.  Since $U$ is a nonempty Scott open set, it is upper and therefore contains $1_M$.  Hence,  the top element $\{1_M\}$ of $H$ belongs to $\BoxU U$.  As $U\neq M$, the bottom element $M$ of $H$ does not belong to $\BoxU U$, so $\BoxU U$ will be proper.
Suppose $K\in\BoxU U$ and $K\leq_H L$.  The reverse-inclusion order means $K\supseteq L$.  Thus $L\subseteq U$ due to the fact that  $K\subseteq U$. So $L\in\BoxU U$.  Consequently,  $\BoxU U$ is upward closed in the frame order.
If $K,L\in\BoxU U$, then
$K\wedge_H L=K\cup L\subseteq U$,
so $K\wedge_H L\in\BoxU U$.  Hence,  $\BoxU U$ is a filter.

It remains to prove Scott openness of  $\BoxU U$ in $H$.   By  \cref{lem:xilawson},  we know $\Sigma M$ is well-filtered.   Since $\BoxU U$ is  open in the upper Vietoris topology,   it is  Scott open by   \cref{lem:smyth-cosar-serscott}.

(ii). According to the above proof, we know that  the upper Vietoris topology is coarser
than $\tau_f$   on $\Q(M)$. By  \cref{lem:xilawson} and   \cref{lem:smyth-cosar-serscott},   the upper Vietoris topology is coarser than the Scott topology on $\Q(M)$. Then the specialization order of $P_S(M)$ and $\Sigma(\Q(M))$  are both the reverse inclusion. Therefore,  the specialization order of $(H,\tau_{\mathrm f})$ is also the reverse inclusion.

(iii).
It is obvious that $\mathord\uparrow x\in \Q(M)=H$ for every $x\in M$, so $\xi$ is well-defined. If $x\leq y$, then
$\mathord\uparrow x\supseteq\mathord\uparrow y$, 
which means
$\xi(x)\leq_{\tau_f}\xi(y)$.
Thus $\xi$ is monotone.

Let $D\subseteq M$ be directed.  Since $M$ is complete, $\bigvee D$ exists, and
$\xi\!\left(\bigvee D\right)=\mathord\uparrow\!\bigvee D=\bigcap\limits_{d\in D}\mathord\uparrow d={\bigvee\limits_{d\in D}}\xi(d)$.
As a result,  $\xi$ preserves directed suprema and is Scott-continuous as a mapping from 
$\Sig M$ to $\Sig H$. As $\tau_{\mathrm f}\subseteq\sigma(H)$,    the   map $\xi$ from $\Sig M$ into the coarser space $(H,\tau_{\mathrm f})$ is continuous.
\end{proof}

The next theorem is the heart of the construction.

\begin{theorem}
\label{frame compact}
Let $M$ be a complete lattice and let
\(H=\Q(M)\)
carry the  topology $\tau_{\mathrm f}$.  If $(H,\tau_{\mathrm f})$ is sober, then the Scott space $\Sig M$ is sober.
\end{theorem}

\begin{proof}
Let $C\subseteq M$ be an arbitrary nonempty irreducible Scott closed subset.  We shall prove that $C$ has a generic point in $\Sig M$.
By \cref{lem:smyth-filter-properties}(iii), the map
\[
   \xi:\Sig M\to(H,\tau_{\mathrm f}),
   \qquad
   x\mapsto\mathord\uparrow x,
\]
is continuous.  A continuous image of an irreducible set is irreducible, and the closure of an irreducible set is irreducible.  Hence, 
$\cl_{\tau_{\mathrm f}}(\xi(C))
$ is a nonempty irreducible closed subset of $(H,\tau_{\mathrm f})$. Let $A=\cl_{\tau_{\mathrm f}}(\xi(C))
$. 
Since  $(H,\tau_{\mathrm f})$ is sober, then there exists $K\in H$ such that
$A=\cl_{\tau_{\mathrm f}}\{K\}$. 
By   \cref{lem:smyth-filter-properties}(ii), we have $\cl_{\tau_{\mathrm f}}\{K\}=\{L\in H:L\supseteq K\}$. 
Therefore,
\begin{equation}\label{eq:A-closure-K}
   A=\{L\in\Q(M):L\supseteq K\}.
\end{equation}

For every $c\in C$,
 $\xi(c)=\mathord\uparrow c\in\xi(C)\subseteq A$. 
By \eqref{eq:A-closure-K},
we can obtain that $\mathord\uparrow c\supseteq K$.
Hence, $k\in\mathord\uparrow c$ holds for every $k\in K$,
which is equivalent to
\begin{equation}\label{eq:C-below-K}
   c\leq k
   \qquad
   \text{for all }c\in C\text{ and }k\in K.
\end{equation}
Thus every element of $K$ is an upper bound of the whole set $C$.

Suppose that $K$ has no least element.  Then for every $k\in K$,  we may choose
$s_k\in K$ such that $k\not\leq s_k$.
For each $k\in K$, it is easy to see that the principal ideal $\mathord\downarrow s_k$ is Scott closed.  Therefore, 
$M\setminus\mathord\downarrow s_k$
is Scott open.
Define
\[
   U
   =
   \bigcup_{k\in K}
      \left(M\setminus\mathord\downarrow s_k\right).
\]
 Moreover, if $k\in K$, our choice $k\not\leq s_k$ gives
$k\in M\setminus\mathord\downarrow s_k\subseteq U$.
Thus
$K\in\BoxU U$.
By \cref{lem:smyth-filter-properties}(i), $\BoxU U$ is a $\tau_{\mathrm f}$-open neighbourhood of $K$.

Since
$K\in A=\cl_{\tau_{\mathrm f}}(\xi(C))$,
every open neighbourhood of $K$ meets $\xi(C)$.  Hence, there exists $c\in C$ such that
$\mathord\uparrow c=\xi(c)\in\BoxU U$.
Equivalently,
$ \mathord\uparrow c\subseteq U$.
In particular, $c\in U$.  Therefore, for some $k\in K$,
$c\in M\setminus\mathord\downarrow s_k$,
so
$ c\not\leq s_k$.
But $s_k\in K$, and \eqref{eq:C-below-K} says $c\leq s_k$.  This contradiction proves that $K$ has a least element.  We denote it by $s$.

As $K$ is saturated in the Scott space, it is an upper set.  It is evident that
$\mathord\uparrow s\subseteq K$.
The reverse inclusion follows because for every $k\in K$, $s\leq k$.  Consequently, 
\begin{equation}\label{eq:K-principal}
   K=\mathord\uparrow s.
\end{equation}

Assume $s\notin C$.  Since $C$ is Scott closed,
$M\setminus C$
is Scott open and contains $s$.  Then  $K=\mathord\uparrow s\subseteq M\setminus C$.
Thus
$K\in\BoxU(M\setminus C)$. By \cref{lem:smyth-filter-properties}(i), $\BoxU(M\setminus C)$ is an open neighbourhood of $K$.  Again using
$K\in\cl_{\tau_{\mathrm f}}(\xi(C))$, there exists $c\in C$ such that
$\mathord\uparrow c\subseteq M\setminus C$.
But $c\in\mathord\uparrow c$ and $c\in C$, a contradiction.  Therefore, 
\begin{equation}\label{eq:s-in-C}
   s\in C.
\end{equation}

As a result,  $\mathord\downarrow s\subseteq C$. By \eqref{eq:C-below-K}, for $c\in C$, $c\leq s$ , so
$C\subseteq\mathord\downarrow s$.
Hence,
$C=\mathord\downarrow s$.
Because  the specialization order of $\Sig M$ is the original order,
$\cl_{\sigma(M)}\{s\}=\mathord\downarrow s$.
Therefore, 
$C=\cl_{\sigma(M)}\{s\}$. It is easy to see that the generic point is unique.  We conclude that every nonempty irreducible Scott closed subset has a unique generic point.  Thus $\Sig M$ is sober.
\end{proof}

\begin{corollary}\label{cor:transfer-contrapositive}
If $M$ is a complete lattice whose Scott space is non-sober, then the  topological space
$(\Q(M),\tau_{\mathrm f})$
is non-sober.
\end{corollary}

\section{The function space with compact-open topology need not   be sober}
\label{sec:main-result}

Let $X$ and $Y$ be topological spaces, and 
 $\mathbb{C}(X,Y)$ be the set of all continuous functions from $X$ to $Y$. The compact-open topology on the set  $\mathbb{C}(X,Y)$  is generated by the subsets of the form $N(K\rightarrow V)=\{f\in\mathbb{C}(X,Y):f(K)\subseteq V\}$, where $K$ is compact in $X$ and $V$ is open in $Y$. Let $[X\rightarrow Y]_K$ denote the function space $\mathbb{C}(X,Y)$ endowed with the compact-open topology. Note that when $Y$ is the Sierpi\'nski space, the function space  $[X\rightarrow Y]_K$ is homeomorphic to the space $(\OO(X),\phi(\OO(X)))$ (see \cite{xubao}).

By the following theorem, we show that  the function space with compact-open topology need  not  be sober.

\begin{theorem}\label{thm:main}
There exists a sober space $X$ such that the space
$(\OO(X),\phi(\OO(X)))$
is not sober.  
\end{theorem}

\begin{proof}
By Isbell's construction \cite{Isbell1982}, there exists a complete lattice $M$ such that
$\Sig M$
is not sober.
Put $H=\Q(M)$,
ordered by the reverse inclusion.  By \cref{Ishey}, $H$ is a spatial frame.
Define
$X=\pt(H)$. For $K\in H$, let
$\widehat K=\{p\in X:p(K)=1\}$.
By \cref{lem:frame-spectrum}, the sets $\widehat K$ form the topology of a sober space $X$ and
\begin{equation}\label{eq:open-frame-H}
   \OO(X)\cong H.
\end{equation}
By \cref{lem:filter-co},  we have 
\begin{equation}\label{eq:function-homeo}
(\OO(X),\phi(\OO(X)))
\cong  (H,\tau_{\mathrm f})=
   (\Q(M),\tau_{\mathrm f}).
\end{equation}

 Suppose  for contradiction that $(\OO(X),\phi(\OO(X)))$ is sober.  By   \eqref{eq:function-homeo}, $(\Q(M),\tau_{\mathrm f})$ is sober.  Then  $\Sig M$ is sober  by \cref{frame compact}.  This contradicts the choice of Isbell's lattice $M$. Therefore $(\OO(X),\phi(\OO(X)))$ is not sober.
\end{proof}

\begin{corollary}
Let $M$ be the Isbell's complete lattice,  $X=\pt(\Q(M))$, and $Y$  the  Sierpi\'nski space. Then  the function space $[X\rightarrow Y]_K$ is not sober.
\end{corollary}

\begin{remark}
The construction does not require an explicit  description of Isbell's lattice.  What is needed is exactly the existence of a complete lattice with non-sober Scott topology.  Once such an $M$ is fixed, the   space $X=\pt(\Q(M))$ is defined canonically from the spatial frame $\Q(M)$.
\end{remark}

\section*{Acknowledgements}

This work is supported by 
the Science and Technology Project of Hebei Education Department (QN2024027), 
the National Natural Science Foundation of China--Mathematics Tianyuan Fund (12526623), and the National Natural Science Foundation of China (12231007, 12301588).

\end{document}